\documentclass[10pt]{article}
\usepackage{amssymb}
\usepackage{mathrsfs}
\usepackage{amsfonts}
\usepackage{tikz}

\usepackage{amsmath,amsthm,amsfonts,amssymb,amscd}

\allowdisplaybreaks[4]
\newtheorem {Problem} {Problem}[section]
\newtheorem {Theorem} [Problem]{Theorem}
\newtheorem {Lemma}[Problem]{Lemma}
\newtheorem{Conjecture}[Problem]{Conjecture}

\newenvironment {Proof}{\noindent {\bf Proof.}}{\hfill\ensuremath{\square}}
\newcommand*{\QEDB}{\hfill\ensuremath{\square}}

\begin{document}

\title{On  the signless Laplacian spectral radius of $K_{s,t}$-minor free graphs  \thanks{This work is supported by  the Joint NSFC-ISF Research Program (jointly funded by the National Natural Science Foundation of China and the Israel Science Foundation (No. 11561141001)),  the National Natural Science Foundation of China (No.11531001) and the Montenegrin-Chinese Science and Technology Cooperation Project (No.3-12).
\newline \indent $^{\dagger}$Corresponding author:
Xiao-Dong Zhang (Email: xiaodong@sjtu.edu.cn),}}

\author{ Ming-Zhu Chen and  Xiao-Dong Zhang
\\
School of Mathematical Sciences, MOE-LSC, SHL-MAC\\
Shanghai Jiao Tong University,
Shanghai 200240, P. R. China\\
Email: chenmingzhuabc@163.com and Xiaodong@sjtu.edu.cn\\
{\bf This paper is  devoted to the memory of Professor Marvin Marcus}
}
\date{}
\maketitle

\begin{abstract}
   In this paper, we prove that if  $G$ is a $K_{2,t}$-minor free graph of order $n\geq t^2+4t+1$ with  $t\geq 3$,  the signless Laplacian spectral radius  $q(G)\leq \frac{1}{2}(n+2t-2+\sqrt{(n-2t+2)^2+8t-8}\ )$ with
equality if and only if $n\equiv 1~(\mathrm{mod}~t)$ and $G=F_{2,t}(n)$, where  $F_{s,t}(n):=K_{s-1}\vee (p\cdot K_t\cup K_r)$  for $n-s+1=pt+r$ and  $0\leq r<t$.  In particular, if $t=3$ and $n\geq 22$,
 then $F_{2,3}(n)$ is the unique $K_{2,3}$-minor free graph of  order $n$ with the maximum  signless Laplacian spectral radius. In addition,  $F_{3,3}(n)$ is the unique  extremal graph  with the maximum  signless Laplacian spectral radius among all $K_{3,3}$-minor free graphs of  order $n\ge 1186$.
\end{abstract}

\section{Introduction}
 Let $G$ be an undirected simple graph with vertex set
$V(G)=\{v_1,\dots,v_n\}$ and edge set $E(G)$.
The \emph{adjacency matrix}
$A(G)$ of $G$  is the $n\times n$ matrix $(a_{ij})$, where
$a_{ij}=1$ if $v_i$ is adjacent to $v_j$, and $0$ otherwise. The  \emph{spectral radius} of $G$ is the largest eigenvalue of $A(G)$. Let  $D(G)$ be the degree diagonal matrix of $G$. The
matrix $Q(G)=D(G)+A(G)$ is known as  the signless Laplacian matrix of $G$. The  \emph{signless Laplacian spectral radius} of $G$ is the largest eigenvalue of $Q(G)$, denoted by
 $q(G)$. For $v\in V(G)$,  the \emph{neighborhood} $N_G(v)$ of $v$  is $\{u: uv\in E(G)\}$ and the \emph{degree} $d_G(v)$ of $v$  is $|N_G(v)|$.
We write $N(v)$ and $d(v)$ for $N_G(v)$ and $d_G(v)$ respectively if there is no ambiguity.
 Denote by $\Delta(G)$ the maximum degree of $G$ and  $\Delta'(G)$ the second largest degree of $G$, i.e., $\Delta(G)=\max\{d(v): v\in V(G)\}$ and   $\Delta'(G)=\max\{d(u): \ u\in V(G)\setminus \{v\}\}$ with $d(v)=\Delta(G)$.  A graph $H$ is called a \emph{minor} of a graph $G$ if it can be obtained from $G$ by deleting edges, contracting edges or deleting vertices. A graph $G$ is \emph{$H$-minor free} if   it  does not contain $H$ as a minor.
 Moreover, a graph $G$ is said to be an {\it edge-maximal $H$-minor free graph} if $G$ is \emph{$H$-minor free} and  the graph $G^+$ which is obtained from $G$ by joining any two nonadjacent vertices of $G$ has a $H$-minor. For two vertex disjoint graphs $G$ and $H$,  we denote by  $G\cup H$ and  $G\vee H$  the \emph{union} of $G$ and $H$,
and the \emph{join} of $G$ and $H$ which is obtained by joining every vertex of $G$ to every vertex of $H$, respectively.
Denote by $k\cdot G$  the  union of $k$ disjoint copies of $G$. For graph notation and terminology undefined here,  readers are referred to \cite{BM}.



The investigation of $H$-minor free graphs is of great significance. It is very useful
for studying the structures and properties of graphs. For example,
Wagner \cite{Wagner} showed that a  graph is planar if and only if it does not contain $K_5$ or  $K_{3,3}$ as  a minor.  Similar to Wagner's theorem, a  graph is outerplanar if and only if it does not contain  $K_4$ or $K_{2,3}$ as a minor \cite{BM}.
In extremal graph theory, one of problems that are concerned about is the the maximum number of edges for graphs that do not contain a given $H$ as a minor. It is known that a planar graph has at most $3n-6$ edges and  an  outerplanar graph has at most $2n-3$ edges, see \cite{BM}. Moreover,  an edge maximal $K_{3,3}$-minor graph has  at most $3n-5$ edges and
 an edge maximal $K_{2,3}$-minor graph has at most $2n-2$ edges, see \cite{Fang, CRS}.
In spectral extremal graph theory, it is interesting to determine the maximum (signless Laplacian) spectral radius of graphs that do not contain a given $H$ as  a minor.
Nikiforov \cite{Nikiforov}  determined all  graphs  maximizing the spectral radius of all  $K_{2,3}$-minor free graphs of large order. In addition,  he established a sharp upper bound  of  spectral radius over all $K_{2,t}$-minor free graphs of large order for $t\geq4$, and determined all extremal graphs.
Recently,  Tait \cite{Tait} extended Nikiforov's result  to  $K_{s,t}$-minor free graphs for $2\leq s\leq t$.
For more details, readers may be referred to \cite{Benediktovich,Hong,SH,Tait,YSH}.

In order to state more results, we need some symbols for given graphs.
Let $F_{s,t}(n):=K_{s-1}\vee (p\cdot K_t\cup K_r)$, where $n-s+1=pt+r$ and $0\leq r<t$. Clearly, $F_{s,t}(n)$ is a $K_{s,t}$-minor free graph of order $n$. 

In \cite{DNP}, it was shown that if $G$ is a $K_{2,2}$-free graph of order $n\geq4$, then  $F_{2,2}(n)$ is the unique $K_{2,2}$-free graph with the maximum  signless Laplacian spectral radius.
Since $F_{2,2}(n)$ is $K_{2,2}$-minor free, $F_{2,2}(n)$ is also the unique $K_{2,2}$-minor free graph with the maximum signless Laplacian spectral radius.

\begin{Theorem}\cite{DNP}\label{K22}
Let $G$ be a $K_{2,2}$-minor free graph of order $n\geq 4$. Then $$q(G)\leq q(F_{2,2}(n))$$ with equality if and only if $G=F_{2,2}(n)$.
\end{Theorem}

Motivated by  above results, we investigate the signless Laplacian spectral radius of  $K_{s,t}$-minor free graphs of order $n$.
For $s=2$ and $ t=3$, we  determine all  graphs  which maximize the signless Laplacian spectral radius  of all $K_{2,3}$-minor free graphs.
 For $s=2$ and $t\geq 4$, we also obtain a sharp upper bound for signless Laplacian spectral radius of $K_{2,t}$-minor free graphs of large order and determine the extremal graphs.
 In addition, for $s=t=3$, we  determine the  graphs of large order $n$ which maximize the signless Laplacian spectral radius  of all $K_{3,3}$-minor free graphs. The main results of this paper are stated as follows.

\begin{Theorem}\label{thm1.1}
Let $G$ be a $K_{2,3}$-minor free graph of order $n\geq 22$. Then $$q(G)\leq q(F_{2,3}(n))$$ with equality if and only if $G=F_{2,3}(n)$, where
$n-1=3p+r$, $0\le r<3$ and $q(F_{2,3}(n))$ is the largest root of the following equation
$$x^3-(n+2r+3)x^2+(8r+3n+2rn-8)x+14r+4n-8rn-2r^2-4=0.$$
\end{Theorem}

\begin{Theorem}\label{thm1.2}
Let $t\geq4$ and $G$ be a $K_{2,t}$-minor free graph of order $n\geq t^2+4t+1$. Then $$q(G)\leq \frac{n+2t-2+\sqrt{(n-2t+2)^2+8t-8}}{2}$$ with
equality if and only if $n\equiv 1~(\mathrm{mod}~t)$ and $G=F_{2,t}(n)$.
\end{Theorem}

\begin{Theorem}\label{thm-K33-minor free}
Let $G$ be a $K_{3,3}$-minor free graph of order $n\geq 1186$. Then $$q(G)\leq q(F_{3,3}(n))$$ with
equality if and only if $G=F_{3,3}(n)$, where $n-2=pt+r$, $0\le r<3$ and $q(F_{3,3}(n))$ is the largest root of the following equation
$$x^3-(n+2r+6)x^2+(12r+4n+2rn+4)x+4r-8rn-4r^2=0.$$

\end{Theorem}

The rest of this paper is organized as follows. In Section~2, we present some known and necessary results.
In Section~3, we give the proofs of Theorems~\ref{thm1.1} and  \ref{thm1.2}.
In Section~4, we present some necessary lemmas and give the proof of Theorem~\ref{thm-K33-minor free}.

\section{Preliminary}


\begin{Lemma}\label{K23-free}\cite{DNP1}
Let $G$ be a $K_{2,3}$-free graph of order $n\geq 22$ with $\Delta(G)\leq n-2$. Then $$q(G)<n.$$
\end{Lemma}

\begin{Lemma}\label{K2t-free}\cite{DNP1}
Let $t\geq2$ and $G$ be a $K_{2,t}$-free graph of order $n\geq t^2+4t+1$. Then $$q(G)\leq \frac{n+2t-2+\sqrt{(n-2t+2)^2+8t-8}}{2}$$ with
equality if and only if $G=K_1\vee H$, where $H$ is a $(t-1)$-regular graph of order $n-1$.
\end{Lemma}

\begin{Lemma}\label{bound for Fst}
(i) For $2\leq s\leq t$ and $n\geq s+2t^2-5$,
$$n+2s-4<q(F_{s,t}(n))\leq \frac{n+2s+2t-6+\sqrt{(n+ 2s - 2t - 2)^2+8(s - 1)(t-s+1)}}{2}.$$
Moreover, the second equality holds if and only if $n\equiv s-1~(\mathrm{mod}~t)$.\\
(ii) For $n-1=3p+r$ and $0\le r<3$, $q(F_{2,3}(n))$ is the largest root of the following equation
$$x^3-(n+2r+3)x^2+(8r+3n+2rn-8)x+14r+4n-8rn-2r^2-4=0.$$
(iii) For $n-2=pt+r$ and $0\le r<3$, $q(F_{3,3}(n))$ is the largest root of the following equation
$$x^3-(n+2r+6)x^2+(12r+4n+2rn+4)x+4r-8rn-4r^2=0.$$
\end{Lemma}

\begin{Proof}
Denote $q=q(F_{s,t}(n))$ and $Q=Q(F_{s,t}(n))$. Let $\mathbf x$ be a positive eigenvector of $Q$ corresponding to $q$.
Let $n-s+1=pt+r$ with $0\le r<t$. By symmetry and the Perron-Frobenius theorem, all vertices of  subgraphs $K_{s-1}$, $p\cdot K_t$, or $K_r$ in  $F_{s,t}(n):=K_{s-1}\vee (p\cdot K_t \cup K_r)$ have the same  eigenvector components respectively, which are denoted by $x_1$, $x_2$, $x_3$, respectively. We consider the following two cases.

{\bf Case 1:}  $r=0$.
By  $Q\mathbf x=q\mathbf x$, it is easy to see that
\begin{eqnarray*}
    (q-n+1) x_1 &=& (s-2)x_1 +(n-s+1)x_2, \\
    (q-s-t+2) x_2 &=& (s-1)x_1+(t-1)x_2.
\end{eqnarray*}
Then $q$ is the largest root of $g(x)=0$, where
$$g(x)=x^2 + (6 - 2s - 2t - n)x + (n + s - 3)(s + 2t - 3) - (s - 1)(n - s + 1).$$
Thus
 $$q= \frac{n+2s+2t-6+\sqrt{(n+ 2s - 2t - 2)^2+8(s-1)(t - s +1)}}{2}.$$

 {\bf Case 2:}  $1\leq r<t$. By  $Q\mathbf x=q\mathbf x$, it is easy to see that
\begin{eqnarray*}
    (q-n+1)x_1 &=& (s-2)x_1 +(n-s-r+1)x_2+rx_3, \\
    (q-s-t+2) x_2 &=& (s-1)x_1+(t-1)x_2,    \\
     (q-s-r+2) x_3 &=& (s-1)x_1+(r-1)x_3.
\end{eqnarray*}
Then $q$ is  the largest root of $f(x)=0$, where
$$f(x)=(x+3-s-2r)g(x)-2r(s - 1)(r - t).$$
Since $q>q(K_{s+r-1})=2s+2r-4$ and $1\leq r<t$, we have
$$g(q)=\frac{2r(s - 1)(r - t)}{q+3-s-2r}<0,$$
which implies that $$q<\frac{n+2s+2t-6+\sqrt{(n+ 2s - 2t - 2)^2+8(s - 1)(t-s+1)}}{2}.$$
Moreover, let $$h(x)=x^2 + (6 - 2s - 2t - n)x + (n + s - 3)(s + 2t - 3) - (s - 1)(n - s).$$
Then $$f(x)=(x+3-s-2r)h(x)-(x+3-s-2r)(s-1)-2r(s - 1)(r - t).$$
Thus $$h(q)=\frac{[q+3-s-2r(1-r+t)](s-1)}{q+3-s-2r}.$$
Note that  $q> q(K_{s-1,n-s+1})=n$, we have $q+3-s-2r>0$ and
$$q+3-s-2r(1-r+t)> n+3-s-2(t-1)(1+t)=n-s-2t^2+5\geq0. $$
Thus $h(q)>0$, which implies that
$q$ is larger than the largest root of $h(x)=0$, i.e.,
$$q> \frac{n+2s+2t-6+\sqrt{(n+ 2s - 2t - 2)^2+4(s - 1)(2t - 2s +1)}}{2}> n+2s-4.$$

 As for (ii) and (iii), the results follow directly from the proof of (i).
\end{Proof}

\begin{Lemma}\label{edge-K1t-minor free}\cite{DJS}
Let $t\geq3$ and  $G$ be a $K_{1,t}$-minor free graph of order $n\geq t+2$. Then $$e(G)\leq n+\frac{1}{2}t(t-3).$$
\end{Lemma}

%

The following Lemmas~\ref{1-k-12}--\ref{largest maximum degree and small second largest degree} are a little different from their original  forms \cite{YWG},  but indeed they  are correct according  to original proofs.

\begin{Lemma}\label{1-k-12}\cite{YWG}
Let $1\leq k\leq12$ and $G$ be a  graph of order $n\geq115$ with degree sequence $d_1\geq \dots\geq d_n$.
If $n-3 \leq d_1\leq n-2$,  $\frac{n}{6}+1\leq d_{k+1}\leq\dots\leq d_2\leq n-61$, and
$d_n\leq \dots \leq d_{k+2}<\frac{n}{6}+1$, then
$$q(G)\leq n+2.$$
\end{Lemma}

\begin{Lemma}\label{small second largest degree}\cite{YWG}
Let  $G$ be a  graph of order $n\geq4$ with degree sequence $d_1\geq \dots\geq d_n$.
If $n-3\leq d_1\leq n-2$ and $d_n\leq \dots \leq d_2<\frac{n}{6}+1$, then
$$q(G)\leq n+2.$$
\end{Lemma}

\begin{Lemma}\label{1-k-13}\cite{YWG}
Let $1\leq k\leq13$ and $G$ be a  graph of order $n\geq91$ with degree sequence $d_1\geq \dots\geq d_n$.
If $d_1=n-1$,  $\frac{n}{7}+\frac{19}{7}\leq d_{k+1}\leq\dots\leq d_2\leq n-75$, and
$d_n\leq \dots \leq d_{k+2}<\frac{n}{7}+\frac{19}{7}$, then
$$q(G)\leq n+2.$$
\end{Lemma}

\begin{Lemma}\label{largest maximum degree and small second largest degree}\cite{YWG}
Let  $G$ be a  graph of order $n\geq6$ with degree sequence $d_1\geq \dots\geq d_n$.
If $ d_1= n-1$ and $d_n\leq \dots \leq d_2<\frac{n}{7}+\frac{19}{7}$, then
$$q(G)\leq n+2.$$
\end{Lemma}

\begin{Lemma}\label{square real matrix}\cite{BP}
Let $A$ be an irreducible nonnegative matrix of order
$n$ and spectral radius $\lambda$. If there exists a nonnegative real vector $y\neq 0$ such that $Ay \leq ry$ $(r \in \mathbb{R})$, then $\lambda\leq r$.
\end{Lemma}

\section{Proofs of Theorems~\ref{thm1.1} and \ref{thm1.2} }

\noindent {\bf Proof of Theorem~\ref{thm1.1}.}
Let $G$ be  a spectral extremal graph with the maximum signless Laplacian spectral radius among all $K_{2,3}$-minor free graphs of order $n\ge 22$.
We just need to prove that $G=F_{2,3}(n)$. 
 Let
$\mathbf{x}=(x_v)_{v\in V(G)}\in \mathbb{R}^n$ be a positive  unit eigenvector corresponding to  $q(G)$.

By Lemma~\ref{bound for Fst}, $q(F_{2,3}(n))>n$ and  $q(F_{2,3}(n))$ is the largest root of  equation
$$x^3-(n+2r+3)x^2+(8r+3n+2rn-8)x+14r+4n-8rn-2r^2-4=0,$$ where $n-1=3p+r$ and $0\le r<3$.
 Since $F_{2,3}(n)$ is $K_{2,3}$-minor free, we have $q(G)\geq q(F_{2,3}(n))>n$.   Hence $\Delta(G)=n-1$ by  Lemma~\ref{K23-free}.  Thus let $u$ be a vertex of $G$ such that $d(u)=\Delta(G)=n-1$.
Clearly every vertex of $G-u$ has degree at most 2 by  $K_{2,3}\nsubseteq G$. In addition, $G-u$ does not contain any cycle of length at least 4 as a subgraph, otherwise $G$ contains $K_{2,3}$ as a minor.  Hence every connected component of  $G-u$  is either a triangle or a path of order at least 1. Furtherer there is at most one connected component of $G-u$ is a path.  Otherwise adding an edge to two pendant vertices in two different connected components which are paths leads to a $K_{2,3}$-minor free graph with larger signless Laplacian spectral radius, a contradiction.
Let $H$ be a connected component of $G-u$ which is not triangle, i.e., $H$ is a path of order $h\geq1$.   Then we have the following claim.

\vspace{2mm}
{\bf Claim:} $1\le h\le 2$.

  Let $v_1,\dots, v_h$ be the vertices along the path. For simplicity, let $x_i=x_{v_i}$.
We consider the following three cases.

\vspace{2mm}
{\bf Case 1:} $h=3$. Then we can add an edge to two endpoints of $H$ and get a $K_{2,3}$-minor free graph with larger signless Laplacian spectral radius, a contradiction.

\vspace{2mm}
{\bf Case 2:}  $h=4$.  By symmetry, $x_1=x_4$ and $x_2=x_3$. Let $G'=G-v_1v_2+v_2v_4$.
Clearly $G'$ is $K_{2,3}$-minor free and
$$q(G')-q(G)\geq(x_2+x_4)^2-(x_1+x_2)^2=0.$$
If $q(G')=q(G)$, then $\mathbf x$ is also an unit eigenvector corresponding to $q(G')$. Since $v_2$, $v_3$, and $v_4$
are  symmetric in $G'$, we have $x_2=x_3=x_4$, which implies that $x_1=x_2=x_3=x_4$.
Considering the eigenequations  of $G$ on vertices $v_1$ and $v_2$, we have
$$ q(G)x_1=2x_1+x_2+x_u \quad \mbox{and} \quad q(G)x_2=3x_2+x_1+x_3+x_u,$$ a contradiction.
Hence $q(G')>q(G)$, a contradiction.

\vspace{2mm}
{\bf Case 3:} $h\geq5$.  If $h$ is odd, say $h=2p+1$ and $p\geq2$.
By symmetry, $x_i=x_{h+1-i}$ for $1\leq i \leq p$. Let $G'=G-\{v_{p-1}v_p,v_{p+2}v_{p+3}\}+\{v_pv_{p+2},v_{p-1}v_{p+3}\}$.
Clearly $G'$ is $K_{2,3}$-minor free and
\begin{eqnarray*}
  &&q(G')-q(G) \\
  &\geq& (x_p+x_{p+2})^2+(x_{p-1}+x_{p+3})^2-(x_{p-1}+x_p)^2-(x_{p+2}+x_{p+3})^2 \\
   &=&  2x_px_{p+2}+ 2x_{p-1}x_{p+3}- 2x_{p-1}x_p- 2x_{p+2}x_{p+3}\\
   &=& 2x_p^2+2x_{p-1}^2-4x_{p-1}x_p\\
   &=&2(x_{p-1}-x_p)^2\\
   &\geq&0.
\end{eqnarray*}
If $h$ is even, say $h=2p$ and $p\geq3$. By symmetry, $x_i=x_{h+1-i}$ for $1\leq i \leq p$. Let $G'=G-\{v_{p-1}v_p,v_{p+2}v_{p+3}\}+\{v_pv_{p+2},v_{p-1}v_{p+3}\}$.
Clearly $G'$ is $K_{2,3}$-minor free and
\begin{eqnarray*}
  &&q(G')-q(G) \\
  &\geq& (x_p+x_{p+2})^2+(x_{p-1}+x_{p+3})^2-(x_{p-1}+x_p)^2-(x_{p+2}+x_{p+3})^2 \\
   &=&  2x_px_{p+2}+ 2x_{p-1}x_{p+3}- 2x_{p-1}x_p- 2x_{p+2}x_{p+3}\\
   &=& 2x_p(x_{p+2}-x_{p-1})-2x_{p+3}(x_{p+2}-x_{p-1})\\
   &=&2(x_p-x_{p+3})(x_{p+2}-x_{p-1})\\
   &=&0.
\end{eqnarray*}
In either case, if $q(G')=q(G)$, then $\mathbf x$ is also an unit eigenvector corresponding to $q(G')$. Since $v_p$, $v_{p+1}$, and $v_{p+2}$ are  symmetric in $G'$, we have $x_{p-1}=x_p=x_{p+1}$. Using the eigenequations of $G$, we have $x_1=\dots=x_h$. Now considering the eigenequations  of $G$ on vertices $v_1$ and $v_2$, we have
$$ q(G)x_1=2x_1+x_2+x_u \quad \mbox{and} \quad q(G)x_2=3x_2+x_1+x_3+x_u,$$ a contradiction.
Hence in either case $q(G')> q(G)$, a contradiction.

Hence Claim holds. Then $G-u$ consists of disjoint copies of triangles and at most a path of order 1 or 2.  So $G=F_{2,3}(n)$.  This completes the proof.
\QEDB

\bigskip
\noindent {\bf Proof of Theorem~\ref{thm1.2}.}
Since $G$ is $K_{2,t}$-minor free, $G$ is $K_{2,t}$-free. By Lemma~\ref{K2t-free},
$$q(G)\leq \frac{n+2t-2+\sqrt{(n-2t+2)^2+8t-8}}{2}$$
and  equality holds if and only if $G=K_1\vee H$, where $H$ is a $(t-1)$-regular graph of order $n-1$.
Hence we just need to prove that the equality in the theorem holds if and only if  $n\equiv 1~(\mathrm{mod}~t)$ and $G=F_{2,t}(n)$,
i.e., $H$ is  the union of disjoint copies of complete graphs of order $t$.

 Suppose that $H$ has a connected component $H_1$ that is not isomorphic to $K_t$ and set $h:=|V(H_1)|$. Clearly $H_1$ is a $(t-1)$-regular graph of order $h\geq t+1$. If $h=t+1$, then any two nonadjacent vertices in $H$ have $t-1$ common neighbours, which combing with the vertex of $K_1$ in $G:=K_1\vee H$ yields  $K_{2,t}$,  a contradiction. Thus $h\geq t+2$. Note that $G$ is $K_{2,t}$-minor free, we have $H_1$ is $K_{1,t}$-minor free.
 By Lemma~\ref{edge-K1t-minor free}, $e(H_1)\leq h+\frac{t(t-3)}{2}$. However, since $H_1$ is  a $(t-1)$-regular graph of order $h$, we have $e(H_1)=\frac{h(t-1)}{2}>h+\frac{t(t-3)}{2}$, a contradiction. Hence $H$ is the union of disjoint copies of  complete graphs of order $t$. This  completes the proof.
\QEDB

\section{Proof of Theorem~\ref{thm-K33-minor free}}
In order to prove Theorem~\ref{thm-K33-minor free}, we first prove the following lemmas.

\begin{Lemma}\label{degree-K33-minor free}
Let $G$ be a $K_{3,3}$-minor free graph of order $n\geq 11$ with $\Delta(G)\leq n-4$. Then $$q(G)\leq n+2.$$
\end{Lemma}

\begin{Proof}
Note that
$$q(G)\leq \max \limits_{v\in V(G)}\bigg\{d(v)+\frac{1}{d(v)}\sum\limits_{w\in N(v)} d(w)\bigg\},$$
which dates back to Merris \cite{ Merris}.
Let $u$ be a vertex of $G$ such that
$$q(G)\leq \max \limits_{v\in V(G)}\bigg\{d(v)+\frac{1}{d(v)}\sum\limits_{w\in N(v)} d(w)\bigg\}=d(u)+\frac{1}{d(u)}\sum\limits_{v\in N(u)}d(v).$$

If $d(u)\leq6$, then $q(G)\leq d(u)+\Delta(G)\leq 6+n-4=n+2$.
So we may  assume that $d(u)\geq7$.
Since $G$ is $K_{3,3}$-minor free,  we have $e(G)\leq3n-5$.
Moreover, $N(u)$ is $K_{2,3}$-minor free and  $e(N(u))\leq2d(u)-2$.
Then
\begin{eqnarray*}
  \sum\limits_{v\in N(u)}d(v) &=& d(u)+2e(N(u))+e\big(N(u),V(G)\backslash (N(u)\cup \{u\})\big)\\
   &\leq& d(u)+2e(N(u))+(e(G)-d(u)-e(N(u))\\
   &=& e(G)+e(N(u))\\
   &\leq& 3n-5+2d(u)-2\\
   &=&3n-7+2d(u).
\end{eqnarray*}
Thus
 $$q(G) \leq d(u)+\frac{1}{d(u)}\sum\limits_{v\in N(u)}d(v)
   \leq d(u)+\frac{3n-7+2d(u)}{d(u)}
   = d(u)+2+\frac{3n-7}{d(u)}.$$
  Let $f(x)= x+2+\frac{3n-7}{x}$, where $7\leq x\leq n-4$.
Since $f(x)$ is convex with respect to $x$, we have
$$ q(G) \leq\max \Big\{f(7),f(n-4)\Big\}
   =\max\bigg\{\frac{3}{7}n+8,n+1+\frac{5}{n-4}\bigg\}
   \leq n+2.$$
This completes the proof.
\end{Proof}

\begin{Lemma}\label{large degree-K33-minor free}
Let $G$ be an edge maximal $K_{3,3}$-minor free graph of order $n\geq 374$ with $n-3\leq\Delta(G)\leq n-2$. Then $$q(G)\leq n+2.$$
\end{Lemma}

\begin{Proof}
If $\Delta'(G)<\frac{n}{6}+1$, then $q(G)\leq n+2$ by Lemma~\ref{small second largest degree}.
So we may assume that that $\Delta'(G)\geq \frac{n}{6}+1$. Next we consider the following two cases.

\vspace{1mm}
{\bf Case~1:} $ \frac{n}{6}+1\leq \Delta'(G)\leq n-61$.  Let $U=\big\{v\in V(G):\frac{n}{6}+1\leq d(v)\leq n-61\big\}$
and $k=|U|$. Clearly $k\geq1$. Since $G$ is an edge maximal $K_{3,3}$-minor free graph, we have $\delta(G)\geq2$ and there is at most one vertex with degree 2, and
$e(G)\leq 3n-5$ by \cite[Theorem~2.4 and Lemma~3.3]{Fang}. Then
\begin{eqnarray*}
  6n-10&\geq&2e(G) = \sum\limits_{v\in U}d(v)+ \sum\limits_{v\in V(G)\setminus U}d(v)\\
& \geq & \Big(\frac{n}{6}+1\Big)k+(n-3)+2+3(n-k-2)\\
   &=& \Big(4+\frac{k}{6}\Big)n-2k-7,
\end{eqnarray*}
which implies that  $1\leq k\leq12$. By Lemma~\ref{1-k-12}, $q(G)\leq n+2$.

\vspace{1mm}
{\bf Case~2:} $\Delta'(G)\geq n-60$. Let $v_1,v_2\in V(G)$ such that $d(v_1)=\Delta(G)$ and $d(v_2)=\Delta'(G)$.
 Let $W=V(G)\backslash (N(v_1)\cup \{v_1\})$. Since $n-3\leq d(v_1)\leq n-2$, there are at most two vertices in $W$.
Let $v$ be any vertex in $V(G)\backslash\{v_1,v_2\}$. We consider the following three subcases.

{\bf Subcase~2.1:} $v, v_2\in N(v_1)$. Then $v$ and $v_2$ have at most two common neighbours in $N(v_1)$, otherwise $K_{3,3}\subseteq G$.
Note that $v$ and $v_2$ can be adjacent to all  vertices in $W\cup \{v_1\}$.
Then $v$ and $v_2$ have at most five common neighbours in $G$. Thus $d(v)+d(v_2)\leq n+5$, which implies that
$$d(v)\leq n+5-d(v_2)\leq n+5-(n-60)=65.$$


{\bf Subcase~2.2:} $v,v_2\in W$.  Then $v$ and $v_2$  have at most two common neighbours in $N(v_1)$, otherwise $K_{3,3}\subseteq G$.
 As a result, $v$ and $v_2$ have at most two common neighbours in $G$. Thus $d(v)+d(v_2)\leq n+2$, which implies that
$$d(v)\leq n+2-d(v_2)\leq n+2-(n-60)=62.$$


{\bf Subcase~2.3:} $v\in N(v_1)$, $v_2\in W$, or $v\in W$, $v_2\in N(v_1)$.
Then $v$ and $v_2$ have at most two common neighbours in $N(v_1)$, otherwise $K_{3,3}\subseteq G$. Note that $v$ and $v_2$ may be adjacent to all  other vertices in $W$.
Then $v$ and $v_2$ have at most three common neighbours in $G$. Thus $d(v)+d(v_2)\leq n+3$, which implies that
$$d(v)\leq n+3-d(v_2)\leq n+3-(n-60)=63.$$

 Hence $d(v)\leq 65$  for any $v\in V(G)\backslash\{v_1,v_2\}$.
Let $\mathbf{x}=(x_v)_{v\in V(G)}\in\mathbb{R}^n$ be a positive vector, where

\[x_v=\left\{
\begin{array}{ccc}
 \vspace{1mm}
  1,&&v\in \{v_1,v_2\}\\
  \frac{ 3}{n-2}, && \mbox{otherwise}.
\end{array}\right.
\]
For any $v\in \{v_1,v_2\}$, we have
$$\frac{d(v)x_v+\sum\limits_{w\in N(v)}x_w}{x_v}\leq  d(v)+1+\frac{3(d(v)-1)}{n-2}\leq n-2+1+\frac{3(n-3)}{n-2}<n+2.$$
For any $v\in V(G)\backslash \{v_1,v_2\}$, we  have
$$\frac{d(v)x_v+\sum\limits_{w\in N(v)}x_w}{x_v}\leq d(v)+\frac{1+1+\frac{3(d(v)-2)}{n-2}}{\frac{3}{n-2}}\leq 65+\frac{2+\frac{3\times63}{n-2}}{\frac{3}{n-2}}\leq n+2.$$
Hence $Q\mathbf{x}\leq(n+2)\mathbf{x}$. By Lemma~\ref{square real matrix}, $q(G)\leq n+2$.
\end{Proof}

\begin{Lemma}\label{LL-degree-K33-minor free}
Let $G$ be an edge maximal $K_{3,3}$-minor free graph of order $n\geq 1186$ with $\Delta(G)=n-1$ and $\Delta'(G)\leq n-2$. Then $$q(G)\leq n+2.$$
\end{Lemma}

\begin{Proof}
If $\Delta'(G)<\frac{n}{7}+\frac{19}{7}$, then $q(G)\leq n+2$  by Lemma~\ref{largest maximum degree and small second largest degree}.
So we may assume  that $\Delta'(G)\geq \frac{n}{7}+\frac{19}{7}$.
Next we consider the following two cases.

\vspace{1mm}
{\bf Case~1:}
$ \frac{n}{7}+\frac{19}{7}\leq \Delta'(G)\leq n-75$. Let $U=\big\{v\in V(G):\frac{n}{7}+\frac{19}{7}\leq d(v)\leq n-75\big\}$
and $k=|U|$. Clearly $k\geq1$. Since $G$ is an edge maximal $K_{3,3}$-minor free graph,  we have $\delta(G)\geq2$ and there is at most one vertex with degree 2, and
$e(G)\leq 3n-5$ by \cite[Theorem~2.4 and Lemma~3.3]{Fang}. Then
\begin{eqnarray*}
  6n-10&\geq&2e(G) =\sum\limits_{v\in U}d(v)+ \sum\limits_{v\in V(G)\backslash U}d(v)\\
& \geq &\Big(\frac{n}{7}+\frac{19}{7}\Big)k+(n-1)+ 2+3(n-k-2)\\
   &=& \Big(4+\frac{k}{7}\Big)n-\frac{2}{7}k-5,
\end{eqnarray*}
which implies that  $1\leq k\leq13$. By Lemma~\ref{1-k-13}, $q(G)\leq n+2$.

\vspace{1mm}
{\bf Case~2:} $n-74\leq\Delta'(G)\leq n-2$. Let $v_1,v_2\in V(G)$ such that $d(v_1)=\Delta(G)$ and $d(v_2)=\Delta'(G)$.
 For any $v\in V(G)\backslash \{v_1,v_2\}$, $v$ and $v_2$ have at most three common neighbours, otherwise $K_{3,3}\subseteq G$. Thus $d(v)+d(v_2)\leq n+3$, which implies that $$d(v)\leq n+3-d(v_2)\leq n+3-(n-74)=77.$$
Let $\mathbf{x}=(x_v)_{v\in V(G)}\in\mathbb{R}^n$ be a positive vector, where

\[x_v=\left\{
\begin{array}{ccc}
  1,&&v=v_1\\
   \vspace{1mm}
   \frac{3}{4},&&v=v_2\\
   \vspace{2mm}
  \frac{ 2}{n-2}, && \mbox{otherwise}.
\end{array}\right.
\]
For simplicity, let $x_i=x_{v_i}$ for $1\leq i\leq2$.  For $v_1$, we have
$$\frac{d(v_1)x_1+\sum\limits_{w\in N(v_1)}x_w}{x_1}\leq  d(v_1)+\frac{3}{4}+\frac{2(d(v_1)-1)}{n-2}= n-1+\frac{3}{4}+\frac{2(n-2)}{n-2}<n+2.$$
For $v_2$, we have
$$\frac{d(v_2)x_2+\sum\limits_{w\in N(v_2)}x_w}{x_2}\leq  d(v_2)+\frac{1+\frac{2(d(v_2)-1)}{n-2}}{\frac{3}{4}}\leq n-2+\frac{1+\frac{2(n-3)}{n-2}}{\frac{3}{4}}<n+2.$$
For any $v\in V(G)\backslash \{v_1,v_2\}$, we  have
$$\frac{d(v)x_v+\sum\limits_{w\in N(v)}x_w}{x_v}\leq d(v)+\frac{1+\frac{3}{4}+\frac{2(d(v)-2)}{n-2}}{\frac{2}{n-2}}\leq 77+\frac{\frac{7}{4}+\frac{2\times75}{n-2}}{\frac{2}{n-2}}\leq n+2.$$
Hence $Q\mathbf{x}\leq(n+2)\mathbf{x}$. By Lemma~\ref{square real matrix}, $q(G)\leq n+2$.
\end{Proof}

\vspace{2mm}

\vspace{3mm}
\noindent{\bf Proof of Theorem~\ref{thm-K33-minor free}.}
Let  $G$ be  a spectral extremal  graph with the maximum  signless Laplacian spectral radius among all $K_{3,3}$-minor free graphs of order $n\ge 1186$.   By the Perron--Fronbenius theorem, $G$ is also an edge maximal $K_{3,3}$-minor free graph. Let
$\mathbf{x}=(x_v)_{v\in V(G)}\in \mathbb{R}^n$ be a positive unit eigenvector corresponding to  $q(G)$.
We just need to prove that $G=F_{3,3}(n)$.

By Lemma~\ref{bound for Fst}, $q(F_{3,3}(n))>n+2$ and $q(F_{3,3}(n))$ is the largest root of the following equation
$$x^3-(n+2r+6)x^2+(12r+4n+2rn+4)x+4r-8rn-4r^2=0,$$ where $n-2=pt+r$ and $0\le r<3$.
 Since $F_{3,3}(n)$ is $K_{3,3}$-minor free,  $q(G)\geq q(F_{3,3}(n))>n+2$.
By Lemmas~\ref{degree-K33-minor free}--\ref{LL-degree-K33-minor free},
$\Delta(G)=\Delta'(G)=n-1$. Hence there exist two vertices $u, v\in V(G)$  such that $d(u)=d(v)=n-1$.
Clearly every vertex of $G-\{u,v\}$ has degree at most 2 by  $K_{3,3}\nsubseteq G$.
In addition, $G-\{u,v\}$ does not contain any cycle of length at least 4 as a subgraph, otherwise $G$ contains $K_{3,3}$ as a minor. Hence  $G-\{u,v\}$ consists of triangles and paths of order at least 1.
Furthermore, since $G$ is edge maximal, there exists at most one connected component in $G-\{u,v\}$  which is not triangle.
Let $H$ be a connected component of $G-\{u,v\}$ which is not triangle, i.e., $H$ is a path of order $h\geq1$.
Then we have the following claim.

\vspace{2mm}
{\bf Claim:} $1\le h\le 2$.

  Let $v_1,\dots, v_h$ be the vertices along the path. For simplicity, let $x_i=x_{v_i}$.
We consider the following three cases.

\vspace{2mm}
{\bf Case 1:} $h=3$. Then we can add an edge to two endpoints of $H$ and get a $K_{3,3}$-minor free graph with larger signless Laplacian spectral radius, a contradiction.

\vspace{2mm}
{\bf Case 2:} $h=4$. By symmetry, $x_1=x_4$ and $x_2=x_3$. Let $G'=G-v_1v_2+v_2v_4$.
Clearly $G'$ is $K_{3,3}$-minor free and
$$q(G')-q(G)\geq(x_2+x_4)^2-(x_1+x_2)^2=0.$$
If $q(G')=q(G)$, then $\mathbf x$ is also a positive unit eigenvector corresponding to  $q(G')$. Since $v_2$, $v_3$, and $v_4$
are  symmetric in $G'$, we have $x_2=x_3=x_4$, which implies that $x_1=x_2=x_3=x_4$.
Considering the eigenequations  of $G$ on vertices $v_1$ and $v_2$, we have
$$ q(G)x_1=3x_1+x_2+x_u+x_v \quad \mbox{and} \quad q(G)x_2=4x_2+x_1+x_3+x_u+x_v,$$ a contradiction.
Hence $q(G')>q(G)$, a contradiction.

\vspace{2mm}
{\bf Case 3:} $h\geq5$.  If $h$ is odd, say $h=2p+1$ and $p\geq2$.
By symmetry, $x_i=x_{h+1-i}$ for $1\leq i \leq p$. Let $G'=G-\{v_{p-1}v_p,v_{p+2}v_{p+3}\}+\{v_pv_{p+2},v_{p-1}v_{p+3}\}$.
Clearly $G'$ is $K_{3,3}$-minor free and
\begin{eqnarray*}
  &&q(G')-q(G) \\
  &\geq& (x_p+x_{p+2})^2+(x_{p-1}+x_{p+3})^2-(x_{p-1}+x_p)^2-(x_{p+2}+x_{p+3})^2 \\
   &=&  2x_px_{p+2}+ 2x_{p-1}x_{p+3}- 2x_{p-1}x_p- 2x_{p+2}x_{p+3}\\
   &=& 2x_p^2+2x_{p-1}^2-4x_{p-1}x_p\\
   &=&2(x_{p-1}-x_p)^2\\
   &\geq&0.
\end{eqnarray*}
If $h$ is even, say $h=2p$ and $p\geq3$. By symmetry, $x_i=x_{h+1-i}$ for $1\leq i \leq p$. Let $G'=G-\{v_{p-1}v_p,v_{p+2}v_{p+3}\}+\{v_pv_{p+2},v_{p-1}v_{p+3}\}$.
Clearly $G'$ is $K_{3,3}$-minor free and
\begin{eqnarray*}
  &&q(G')-q(G) \\
  &\geq& (x_p+x_{p+2})^2+(x_{p-1}+x_{p+3})^2-(x_{p-1}+x_p)^2-(x_{p+2}+x_{p+3})^2 \\
   &=&  2x_px_{p+2}+ 2x_{p-1}x_{p+3}- 2x_{p-1}x_p- 2x_{p+2}x_{p+3}\\
   &=& 2x_p(x_{p+2}-x_{p-1})-2x_{p+3}(x_{p+2}-x_{p-1})\\
   &=&2(x_p-x_{p+3})(x_{p+2}-x_{p-1})\\
   &=&0.
\end{eqnarray*}
In either case, if $q(G')=q(G)$, then $\mathbf x$ is also a positive unit eigenvector corresponding to  $q(G')$. Since $v_p$, $v_{p+1}$, and $v_{p+2}$ are  symmetric in $G'$, we have $x_{p-1}=x_p=x_{p+1}$. Using the eigenequations of $G$, we have $x_1=\dots=x_h$. Now considering the eigenequations  of $G$ on vertices $v_1$ and $v_2$, we have
$$ q(G)x_1=3x_1+x_2+x_u+x_v \quad \mbox{and} \quad q(G)x_2=4x_2+x_1+x_3+x_u+x_v ,$$ a contradiction.
Hence in either case $q(G')> q(G)$, a contradiction.

Hence Claim holds. Then $G-\{u,v\}$ consists of disjoint copies of  triangles and at most a path of order 1 or 2.  So $G=F_{3,3}(n)$. This completes the proof.
\QEDB

\vspace{3mm}
It is interesting to see that  if $G$ is a $K_{2,t}$-minor free graph of order $n$, the spectral extremal graphs with the maximum spectral radius and the maximum signless Laplacian spectral radius coincide.  In \cite{Tait}, Tait proposed a general conjecture of spectral radius for $K_{s,t}$-minor free graph of large order $n$.
So we may propose a similar conjecture of signless Laplacian spectral radius for $K_{s,t}$-minor free graph of large order $n$.

\begin{Conjecture}
Let $2\leq s\leq t$ and $G$ be a $K_{s,t}$-minor free graph of sufficiently large order $n$. Then $$q(G)\leq q(F_{s,t}(n))$$ with
equality if and only if $G=F_{s,t}(n)$.
\end{Conjecture}

\indent{\bf Acknowledgements:}

  The authors would like to thank the anonymous referee for many helpful and constructive suggestions to an earlier version of this paper, which results in a great improvement.

\end{document}